\documentclass{article}
\usepackage{array}
\usepackage{amsmath}
\usepackage{amssymb}
\usepackage{graphicx}
\usepackage{url}
\usepackage{bbm}
\usepackage{multicol} 
\usepackage{maplestd2e}

\DefineParaStyle{Maple Heading 1}
\DefineParaStyle{Maple Text Output}
\DefineParaStyle{Maple Dash Item}
\DefineParaStyle{Maple Bullet Item}
\DefineParaStyle{Maple Normal}
\DefineParaStyle{Maple Heading 4}
\DefineParaStyle{Maple Heading 3}
\DefineParaStyle{Maple Heading 2}
\DefineParaStyle{Maple Warning}
\DefineParaStyle{Maple Title}
\DefineParaStyle{Maple Error}
\DefineCharStyle{Maple Hyperlink}
\DefineCharStyle{Maple 2D Math}
\DefineCharStyle{Maple Maple Input}
\DefineCharStyle{Maple 2D Output}
\DefineCharStyle{Maple 2D Input}
\usepackage[top=.9in,bottom=.9in]{geometry}
 \numberwithin{equation}{section}
 
\DefineParaStyle{Maple Output}
\DefineCharStyle{2D Comment}
\DefineCharStyle{2D Math}
\DefineCharStyle{2D Output}
\begin{document}
\pagestyle{plain}
\numberwithin{equation}{section}
\begin{flushleft}

\centerline{\bf On Quadratic Gauss Sums and Variations Thereof} \vskip .2in

\centerline{ M.L. Glasser}\vskip .2in

\centerline{Department of Physics, Clarkson University}
 
\centerline{Potsdam, NY 13699-5820}

\centerline{Dpto. de F\'isica Te\'orica, Universidad de Valladolid}
\centerline{Valladolid 47011, Spain}
\centerline{Email: laryg@clarkson.edu}
\centerline{}
\centerline{and}\vskip .2in
\centerline{ Michael Milgram}

\centerline{Consulting Physicist, Geometrics Unlimited, Ltd.}
\centerline{Box 1484, Deep River, Ont. Canada. K0J 1P0}
\centerline{Email: mike@geometrics-unlimited.com}
\centerline{}
\centerline{Oct. 22, 2014}
\centerline{}

\vskip .2in
{\bf ABSTRACT}\vskip .1in

A number of new terminating series involving $\sin(n^2/k)$ and $\cos(n^2/k)$ are presented and connected to Gauss quadratic sums. Several new closed forms of generic Gauss quadratic sums are obtained and previously known results are generalized.
\newline
\section{Introduction}

In a recent work \cite{GM}, a number of new integrals were evaluated analytically, and in the process, we noted that the limiting case of some of those integrals reduced to simpler known forms that involved trigonometric series with quadratic dummy indices of summation. A reasonably thorough search of the literature indicates that such series are not very well tabulated - in fact only two were found in tables \cite{Hansen} - those two corresponding to the real and imaginary part of Gauss quadratic sums, the series associated with classical number theory \cite{Apostol}, and, recently applied to quantum mechanics \cite{ArmRog},\cite{GhLoo}. Hardy and Littlewood \cite{HardyLitt}, in their original work considered integrals that are similar to those considered here (see \eqref{X6A} below) but not in closed form; a comprehensive historical summary can be found in \cite[Section 2.3]{BerEvans}.\newline

In Section 2, we summarize the integrals upon which our results are based; all eight variations of alternating quadratic sums analogous to  classical quadratic Gauss sums are developed in Section 3.1. These new sums are then connected to the classical sums (Section 3.2), and that relationship is used to find new closed forms for odd-indexed variants of the classical sums (Section 3.3), all of which are believed to be new. The upper summation limit of an interesting subset of these results is extended by the addition of a new parameter $p$ (Section 4.1), and a template proof is provided in the Appendix. Miscellaneous results arising during our investigation are also listed (Section 4.2). The trigonometric form of the results presented in the first four sections are finally rewritten in the canonical form of a generalized Gauss quadratic sum involving the roots of unity, again leading to new expressions, best described as alternating and extended Gauss quadratic sums (Section 5). It is recognized that many of the results derived here can be obtained by clever (and lengthy) manipulation \cite{Chapman} of known results from number theory starting from the classical Gauss sum. This work provides a simple method of obtaining sums that are not listed in any of the tables or reference works usually consulted (e.g. \cite{Hansen}), particularly motivated by the fact that number theoretic methods are not well-known in other fields where such sums arise \cite{ArmRog},\cite{GhLoo}.   

\section{The Basic Integrals}

In the previous work \cite[Eq.(3.33)]{GM}, with $k\in \mathbb{N}, \{ a,b,s\} \in \mathfrak{R}$ we evaluated a family of related integrals of the form

\begin{equation}
\mapleinline{inert}{2d}{X6A :=
Int(v^s*(v^2+1)^(1/2*s)*(sin(4*a*v^2)*cos(s*arctan(1/v))*cosh(
4*a*v)-cos(4*a*v^2)*sin(s*arctan(1/v))*sinh(4*a*v))/cosh(Pi*(2*k-1)*v)
,v = 0 .. infinity) =
-(-1)^k*Sum(((2*k-2*n-1)*(2*k+2*n-1)/(2*k-1)^2)^s*(-1)^n*sin(a*(2*k-2*
n-1)*(2*k+2*n-1)/(2*k-1)^2),n = 1 ..
k-1)/(2*k-1)-1/2*(-1)^k*sin(a)/(2*k-1);}{%
\maplemultiline{
{\displaystyle \int _{0}^{\infty }} 
{\displaystyle \frac {v^{s}\,(v^{2} + 1)^{s/2}\,(
\mathrm{sin}(a\,v^{2})\,\mathrm{cos}(s\,\mathrm{arctan}(
{\displaystyle 1/v} ))\,\mathrm{cosh}(a\,v) - 
\mathrm{cos}(a\,v^{2})\,\mathrm{sin}(s\,\mathrm{arctan}(
{\displaystyle 1/v} ))\,\mathrm{sinh}(a\,v))}{
\mathrm{cosh}(\pi \,(2\,k - 1)\,v)}} dv 
\\ = {\displaystyle {\frac{(-1)^{k+1}}{2^{2\,s}\,(2\,k - 1)^{(2\,s+1)}} \, {\displaystyle 
\sum _{n=1}^{k - 1}} \,(-1)^{n}\,{\displaystyle {((2\,k-1)^2 - (2\,n)^2)\,
}{}}^{ {s}}\,\mathrm{sin}
({\displaystyle \frac {a\,((2\,k-1)^2 - (2\,n)^2)}{4\,(
2\,k - 1)^{2}}} ) }{}}  
\mbox{}- {\displaystyle \frac {1}{2}} \,{\displaystyle \frac {(
-1)^{k}\,\mathrm{sin}(a/4)}{2\,k - 1}}  }
}
\label{X6A}
\end{equation}
where $4a \leq (2 k-1) \pi$, and, for particular choices of $s$, this family of integrals \cite[Eqs. (3.28), (3.34), (3.46 with s=0 and $b\rightarrow i\,a$) and (3.37 and 3.41 combined)]{GM} gives rise (respectively) to the following particular results :
\begin{equation}
\mapleinline{inert}{2d}{R2 := Int(cos(a*v^2)*sinh(2*a*v)/sinh(k*Pi*v),v = 0 .. infinity) =
(-1)^(k+1)*Sum((-1)^n*sin(a*(k+n)*(k-n)/k^2),n = 1 ..
k-1)/k-1/2*(-1)^k*sin(a)/k;}{%
\maplemultiline{
{\displaystyle \int _{0}^{\infty }} 
{\displaystyle \frac {\mathrm{cos}(a\,v^{2})\,\mathrm{sinh}(2\,a
\,v)}{\mathrm{sinh}(k\,\pi \,v)}} \,dv= 
 {\displaystyle \frac {(-1)^{k+1}}{k} \, {\displaystyle 
\sum _{n=1}^{k - 1}} \,(-1)^{n}\,\mathrm{sin}({\displaystyle 
\frac {a\,(k^2 - n^2)}{k^{2}}} )  }  + 
{\displaystyle } \,{\displaystyle \frac {(-1)^{k+1}\,
\mathrm{sin}(a)}{2\,k}}  }
\label{A1R1}
}
\end{equation}

\begin{equation}
\mapleinline{inert}{2d}{Int(exp(i*a*x^2)*cosh(a*x)/cosh((2*k-1)*Pi*x),x = 0 .. infinity) =
(-1)^(k+1)*exp(1/4*i*a)*(1/2+Sum((-1)^n*exp(-i*a*n^2/(2*k-1)^2),n = 1
.. k-1))/(2*k-1);}{%
\[
{\displaystyle \int _{0}^{\infty }} {\displaystyle \frac {e^{i\,
a\,v^{2}}\,\mathrm{cosh}(a\,v)}{\mathrm{cosh}((2\,k - 1)\,\pi \,
v)}} \,dv={\displaystyle \frac{(-1)^{(k + 1)}\,\exp({{i\,a/4}{
})}}{2\,k - 1}{\displaystyle \left(  \! {\displaystyle \frac {1}{2}}  +  
{\displaystyle \sum _{n=1}^{k - 1}} \,(-1)^{n}\,\exp{( - \frac {i\,
a\,n^{2}}{(2\,k - 1)^{2}})}   \!  \right) }{}
} 
\]
}
\label{XY6B}
\end{equation}
\begin{equation}
\mapleinline{inert}{2d}{Int(exp(i*a*v^2)*sinh(a*v)/sinh(2*Pi*k*v),v = 0 .. infinity) =
1/4*i*(-1)^(k+1)*sum((-1)^n*exp(1/4*i*a*(k+1+n)*(k-1-n)/k^2),n = -k ..
-1+k)/k;}{%
\[
{\displaystyle \int _{0}^{\infty }} {\displaystyle \frac {e^{i\,
a\,v^{2}}\,\mathrm{sinh}(a\,v)}{\mathrm{sinh}(2\,\pi \,k\,v)}} 
\,dv={\displaystyle \frac {i}{4\,k}} \,{\displaystyle  {(-1)
^{k }\,\exp(i\,a/4)  {\displaystyle \sum _{n=1 - k}^{ k }
} \,(-1)^{n}\,\exp{(\frac {-i\,a\,n^{2}}{4\,k^{2}
})} }{}} 
\]
}
\label{C10b}
\end{equation}

\begin{equation}
\mapleinline{inert}{2d}{C7A2a :=
Int(cosh(Pi*(2*k-1)*v)*cos(a*v^2)*cosh(a*v)/(cosh(2*Pi*b)+cosh(4*Pi*v*
k-2*Pi*v)),v = 0 .. infinity) =
-1/4*Sum((-1)^n*cosh(2*a*b*(-n+k)/(2*k-1)^2)*cos(1/4*a*(4*b^2-4*n^2+1+
8*k*n-4*k)/(2*k-1)^2),n = 1 .. 2*k-1)/(2*k-1)/cosh(Pi*b);}{%
\maplemultiline{
{\displaystyle \int _{0}^{\infty }} 
{\displaystyle \frac {\mathit{e}^{i\,a\,v^{2}}\,\mathrm{cosh}(\pi \,(2\,k - 1)\,v)\,
\mathrm{cosh}(a\,v)}{\mathrm{cosh}(2\,\pi
 \,b) + \mathrm{cosh}(2\,\pi\,v\,(2\,k-1))}} \,dv= 
 - {\displaystyle \frac {1}{4\,(2\,k - 1)\,\mathrm{cosh}
(\pi \,b)}} \,{\displaystyle {
{\displaystyle \sum _{n=1}^{2\,k - 1}} \,(-1)^{n}\,\mathrm{cosh}(X_{n} )\,
\mathit{e}^{i\,A_{n}}}\,{}}  }
}
\label{C7A2}
\end{equation}
where
\begin{equation}
\mapleinline{inert}{2d}{A[n] = 1/4*a*(4*b^2-4*n^2+1+8*k*n-4*k)/(2*k-1)^2;}{%
\[
{A_{n}}={\displaystyle \frac {a\,(4\,b^{2} - 4\,n^{2} + 1 + 8\,k
\,n - 4\,k)}{4\,(2\,k - 1)^{2}}} 
\]
}
\end{equation}
and
\begin{equation}
\mapleinline{inert}{2d}{X[n] = 2*a*b*(-n+k)/(2*k-1)^2;}{%
\[
{X_{n}}={\displaystyle \frac {2\,a\,b\,( - n + k)}{(2\,k - 1)^{2}
}.} 
\]
}
\end{equation}

\section{Basic Derivations}
\subsection{Alternating Gauss Quadratic Sums}
All the main results in this section depend upon the following well-known result \cite[Eq.3.691]{G&R}

\begin{equation}
\mapleinline{inert}{2d}{Int(sin(a*v^2), v = 0 .. infinity) = (1/4)*sqrt(2)*sqrt(Pi)/sqrt(a)}{\[\displaystyle \int _{0}^{\infty }\!\sin \left( a{v}^{2} \right) {dv}=\int _{0}^{\infty }\!\cos \left( a{v}^{2} \right) {dv}=1/4\,\sqrt{\frac { {2\,\pi }}{ {a}}}\]}\,,\hspace{2cm} a>0.
\end{equation}
In \eqref{A1R1}, we evaluate the limit $a=k\pi/2$, and, after setting $k \rightarrow 2 k$, we eventually arrive at: 

\begin{equation}
\mapleinline{inert}{2d}{R2E := sum(sin(1/4*Pi*n^2/k)*(-1)^n,n = 1 .. 2*k-1) =
(-1)^k*sqrt(1/2*k);}{%
\[
{\displaystyle \sum _{n=1}^{2\,k }} \,(-1)^{n} \,\mathrm{
sin}({\displaystyle \frac {\pi \,n^{2}}{4\,k}} )=(-1)^{
k}\,\sqrt{{\displaystyle \frac {k}{2}} }\, ;
\]
}
\label{A1R2e}
\end{equation}

after setting $k \rightarrow 2 k-1$, we find

\begin{equation}
\mapleinline{inert}{2d}{sum((-1)^n*cos(1/2*Pi*n^2/(2*k-1)),n = 1 .. 2*k-2) =
-1/2*(-1)^k*(2*k-1)^(1/2)-1/2;}{%
\[
{\displaystyle \sum _{n=1}^{2\,k - 2}} \,(-1)^{n}\,\mathrm{cos}(
{\displaystyle \frac {\pi \,n^{2}}{2\,(2\,k - 1)}} )= - 
{\displaystyle \frac {(-1)^{k}\,\sqrt{2\,k - 1}}{2}}  - 
{\displaystyle \frac {1}{2}} \,.
\]
}
\label{A1R2o}
\end{equation}

Both of these results may be characterized as ${Alternating \;Gauss\; Quadratic\; Sums}$ for which we can find no references in the literature (e.g. \cite{Apostol}). In \eqref{XY6B} and \eqref{C10b} let $a=(2k-1)\pi$ and $a=2k\pi$ respectively, and, by comparing real and imaginary parts, eventually arrive at 
\begin{equation}
\mapleinline{inert}{2d}{sum((-1)^n*cos(1/4*Pi*n^2/k),n = 1 .. 2*k) =
1/H*(-1+(-1)^k+sqrt(2)*sqrt(k)*(-1)^k);}{%
\[
{\displaystyle \sum _{n=1}^{2\,k}} \,(-1)^{n}\,\mathrm{cos}(
{\displaystyle \frac {\pi \,n^{2}}{4\,k}} )={\displaystyle 
\frac { - 1 + (-1)^{k} + \sqrt{2}\,\sqrt{k}\,(-1)^{k}}{2}} 
\]
}
\label{Ans1b}
\end{equation}
and
\begin{equation}
\mapleinline{inert}{2d}{Sum((-1)^n*sin(1/2*Pi*n^2/(2*k-1)),n = 1 .. 2*k-1) =
1/2*(-1)^k*(1+(2*k-1)^(1/2));}{%
\[
{\displaystyle \sum _{n=1}^{2\,k - 1}} \,(-1)^{n}\,\mathrm{sin}(
{\displaystyle \frac {\pi \,n^{2}}{2\,(2\,k - 1)}} )=
{\displaystyle \frac {(-1)^{k}\,(1 + \sqrt{2\,k - 1})}{2}} \; . 
\]
}
\label{Ans3b}
\end{equation}

All of these are believed to be new.\newline

By setting $s=0, a=(2k-1) \pi/4$ in \eqref{X6A} we obtain companions to the above:

\begin{equation} 
\mapleinline{inert}{2d}{sum((-1)^n*cos(Pi/(4*k-1)*n^2),n = 1 ..
2*k-1)+sum((-1)^n*sin(Pi/(4*k-1)*n^2),n = 1 .. 2*k-1) =
-1/2+1/2*(-1)^k*(4*k-1)^(1/2);}{%
\[
 {\displaystyle \sum _{n=1}^{2\,k - 1}} \,(-1)^{n}\,
\left( \! \mathrm{cos}({\displaystyle \frac {\pi \,n^{2}}{4\,k - 1}} )  
 +   \mathrm{sin}({\displaystyle \frac {\pi \,n^{2}}{4\,k - 
1}} )\! \right)  = - {\displaystyle \frac {1}{2}}  + 
{\displaystyle \frac {1}{2}} \,(-1)^{k}\,\sqrt{4\,k - 1}
\]
}
\label{X6Ca}
\end{equation}
after setting $k \rightarrow 2k$, and 
\begin{equation}
\mapleinline{inert}{2d}{sum((-1)^n*cos(Pi/(4*k-3)*n^2),n = 1 ..
2*k-2)-sum((-1)^n*sin(Pi/(4*k-3)*n^2),n = 1 .. 2*k-2) =
-1/2-1/2*(-1)^k*(4*k-3)^(1/2);}{%
\[
 {\displaystyle \sum _{n=1}^{2\,k - 2}} \,(-1)^{n}\,
\left( \mathrm{cos}({\displaystyle \frac {\pi \,n^{2}}{4\,k - 3}} ) 
 -  \mathrm{sin}({\displaystyle \frac {\pi \,n^{2}}{4\,k - 
3}} )\right )  = - {\displaystyle \frac {1}{2}}  - 
{\displaystyle \frac {1}{2}} \,(-1)^{k}\,\sqrt{4\,k - 3}
\]
}
\label{X6Cb}
\end{equation}
after setting $k \rightarrow 2 k -1$.
\newline

If $a=(2k-1) \pi$ in \eqref{XY6B}, after setting $k \rightarrow 2k$ we find
\begin{equation}
\mapleinline{inert}{2d}{Sum((-1)^n*(cos(Pi*n^2/(4*k-1))-sin(Pi*n^2/(4*k-1))),n = 1 .. 2*k-1)
= -1/2*(-1)^k*(4*k-1)^(1/2)-1/2;}{%
\[
{\displaystyle \sum _{n=1}^{2\,k - 1}} \,(-1)^{n}\,\left( \mathrm{cos}(
{\displaystyle \frac {\pi \,n^{2}}{4\,k - 1}} ) - \mathrm{sin}(
{\displaystyle \frac {\pi \,n^{2}}{4\,k - 1}} ) \right ) = - 
{\displaystyle \frac {1}{2}} - 
{\displaystyle \frac {(-1)^{k}\,\sqrt{4\,k - 1}}{2}}  
\]
}
\label{Y6Ca2}
\end{equation}

and, after setting $k \rightarrow 2k-1$ we obtain
\begin{equation}
\mapleinline{inert}{2d}{Y6Cb2 := Sum((-1)^n*(cos(Pi*n^2/(4*k-3))+sin(Pi*n^2/(4*k-3))),n = 1
.. 2*k-2) = -1/2-1/2*(-1)^k*(4*k-3)^(1/2);}{%
\[
{\displaystyle \sum _{n=1}^{2\,k - 2}} \,(-1)^{
n}\,\left ( \mathrm{cos}({\displaystyle \frac {\pi \,n^{2}}{4\,k - 3}} )
 + \mathrm{sin}({\displaystyle \frac {\pi \,n^{2}}{4\,k - 3}} ) \right ) =
 - {\displaystyle \frac {1}{2}}  - {\displaystyle \frac {(-1)^{k}
\,\sqrt{4\,k - 3}}{2}}\,\,. 
\]
}
\label{Y6Cb2}
\end{equation}

Adding and subtracting \eqref{X6Ca} and  \eqref{Y6Ca2} gives

\begin{equation}
\mapleinline{inert}{2d}{Sum((-1)^n*cos(Pi/(4*k-1)*n^2),n = 1 .. 2*k-1) = -1/2;}{%
\[
{\displaystyle \sum _{n=1}^{2\,k - 1}} \,(-1)^{n}\,\mathrm{cos}(
{\displaystyle \frac {\pi \,n^{2}}{4\,k - 1}} )={\displaystyle 
-1/2} 
\]
}
\label{S1}
\end{equation}
and

\begin{equation}
\mapleinline{inert}{2d}{Sum((-1)^n*sin(Pi/(4*k-1)*n^2),n = 1 .. 2*k-1) =
1/2*(-1)^k*(4*k-1)^(1/2);}{%
\[
{\displaystyle \sum _{n=1}^{2\,k - 1}} \,(-1)^{n}\,\mathrm{sin}(
{\displaystyle \frac {\pi \,n^{2}}{4\,k - 1}} )={\displaystyle 
\frac {1}{2}} \,(-1)^{k}\,\sqrt{4\,k - 1}\,\,,
\]
}
\label{S2}
\end{equation}

while performing the same operations on \eqref{X6Cb} and \eqref{Y6Cb2} yields
\begin{equation}
\mapleinline{inert}{2d}{Sum((-1)^n*cos(Pi*n^2/(4*k-3)),n = 1 .. 2*k-2) =
-1/2-1/2*(-1)^k*(4*k-3)^(1/2);}{%
\[
{\displaystyle \sum _{n=1}^{2\,k - 2}} \,(-1)^{n}\,\mathrm{cos}(
{\displaystyle \frac {\pi \,n^{2}}{4\,k - 3}} )= - 
{\displaystyle \frac {1}{2}}  - {\displaystyle \frac {(-1)^{k}\,
\sqrt{4\,k - 3}}{2}} 
\]
\label{S3}
}
\
\end{equation}
and
\begin{equation}
\mapleinline{inert}{2d}{Sum((-1)^n*sin(Pi*n^2/(4*k-3)),n = 1 .. 2*k-2) = 0;}{%
\[
{\displaystyle \sum _{n=1}^{2\,k - 2}} \,(-1)^{n}\,\mathrm{sin}(
{\displaystyle \frac {\pi \,n^{2}}{4\,k - 3}} )=0
\]
}
\label{S4}
\end{equation}
all of which, together with \eqref{A1R2e}, \eqref{A1R2o}, \eqref{Ans1b} and \eqref{Ans3b} are believed to be new; note the fourfold modularity of these results, and the fact that \eqref{S1}, \eqref{S2},\eqref{S3} and \eqref{S4} can be also written, respectively, (and more conveniently)
\begin{equation}
\mapleinline{inert}{2d}{Sum((-1)^n*cos(Pi/(4*k-1)*n^2),n = 1 .. 2*k-1) = -1/2;}{%
\[
{\displaystyle \sum _{n=1}^{4\,k - 1}} \,(-1)^{n}\,\mathrm{cos}(
{\displaystyle \frac {\pi \,n^{2}}{4\,k - 1}} )={\displaystyle 
0}\, 
\]
}
\label{C1a}
\end{equation}

\begin{equation}
\mapleinline{inert}{2d}{Sum((-1)^n*cos(Pi/(4*k-1)*n^2),n = 1 .. 2*k-1) = -1/2;}{%
\[
{\displaystyle \sum _{n=1}^{4\,k - 1}} \,(-1)^{n}\,\mathrm{sin}(
{\displaystyle \frac {\pi \,n^{2}}{4\,k - 1}} )={\displaystyle 
} \,(-1)^{k}\,\sqrt{4\,k - 1}\,\,
\]
}
\label{S1a}
\end{equation}

\begin{equation}
\mapleinline{inert}{2d}{Sum((-1)^n*cos(Pi*n^2/(4*k-3)),n = 1 .. 4*k-3) =
-(-1)^k*(4*k-3)^(1/2);}{%
\[
{\displaystyle \sum _{n=1}^{4\,k - 3}} \,(-1)^{n}\,\mathrm{cos}(
{\displaystyle \frac {\pi \,n^{2}}{4\,k - 3}} )= - (-1)^{k}\,
\sqrt{4\,k - 3}
\]
}
\label{S3a}
\end{equation}
\begin{equation}
\mapleinline{inert}{2d}{Sum((-1)^n*sin(Pi*n^2/(4*k-3)),n = 1 .. 2*k-2) = 0;}{%
\[
{\displaystyle \sum _{n=1}^{4\,k - 3}} \,(-1)^{n}\,\mathrm{sin}(
{\displaystyle \frac {\pi \,n^{2}}{4\,k - 3}} )=0
\]
}
\label{S4a}
\end{equation}  
\newline
%
%

\subsection{Variations}

An interesting finite sum identity/equality arises from \eqref{C7A2}. Set $b=0,\,    a=(2k-1)\,\pi$ in the real and imaginary parts to  discover
\begin{equation}
\mapleinline{inert}{2d}{Sum((-1)^n*cos(1/4*Pi*(2*n-1)*(-2*n-1+4*k)/(2*k-1)),n = 1 .. 2*k-1) =
-sqrt(k-1/2);}{%
\[
{\displaystyle \sum _{n=1}^{2\,k - 1}} \,(-1)^{n}\,\mathrm{sin}(
{\displaystyle \frac {\pi \,(2\,n - 1)\,( - 2\,n - 1 + 4\,k)}{4\,
(2\,k - 1)}} )={\displaystyle \sum _{n=1}^{2\,k - 1}} \,(-1)^{n}\,\mathrm{cos}(
{\displaystyle \frac {\pi \,(2\,n - 1)\,( - 2\,n - 1 + 4\,k)}{4\,
(2\,k - 1)}} )= - \sqrt{k - {\displaystyle \frac {1}{2}} }\,.
\]
}
\label{S1=C1}
\end{equation}
By removing odd multiples of $\pi/2$ from the arguments of each of the above trigonometric functions, \eqref{S1=C1} can alternatively be written in the form of \it{odd-indexed quadratic Gauss sums} \rm
\begin{equation}
\mapleinline{inert}{2d}{sum(sin(1/4*Pi*(2*n-1)^2/(2*k-1)),n = 1 .. 2*k-1) = sqrt(k-1/2);}{%
\[
{\displaystyle \sum _{n=1}^{2\,k - 1}} \,\mathrm{sin}(
{\displaystyle \frac {1}{4}} \,{\displaystyle \frac {\pi \,(2\,n
 - 1)^{2}}{2\,k - 1}} )=
\]
}
\mapleinline{inert}{2d}{sum(cos(1/4*Pi*(2*n-1)^2/(2*k-1)),n = 1 .. 2*k-1) = sqrt(k-1/2);}{%
\[
{\displaystyle \sum _{n=1}^{2\,k - 1}} \,\mathrm{cos}(
{\displaystyle \frac {1}{4}} \,{\displaystyle \frac {\pi \,(2\,n
 - 1)^{2}}{2\,k - 1}} )=\sqrt{k - {\displaystyle \frac {1}{2}} }\,.
\]
\label{S1a=C1a}
}
\end{equation}

The fact that the arguments of both the sine and cosine terms are fractional multiples of $\pi/4$ may explain the intriguing equality of these sine and cosine sums of equal argument, which also hints at the existence of some more general property. The only previously tabulated sums of this genre, corresponding to Gauss' classical result for quadratic sums \cite[Eq. 9.10(30)] {Apostol}, \cite{BerEvans} are to be found in the following tabular listing \cite[Eqs. (16.1.1) and (19.1.1)] {Hansen}: 
\begin{equation}
\sum_{n=1}^{k} \sin(\frac{2\pi n^{2}}{k}) = \frac{\sqrt{k}} {2}\left( 1 + \cos(k\pi/2)-\sin(k\pi /2) \right),
\label{Hs1}
\end{equation}
and
\begin{equation}
\sum_{n=1}^{k} \cos(\frac{2\pi n^{2}}{k}) = \frac{\sqrt{k}} {2}\left( 1 + \cos(k\pi/2)+\sin(k\pi /2) \right).
\label{Hc1}
\end{equation}
This will be discussed further in section 5.
\newline


By adding and subtracting the above, and applying simple trigonometric identities, we also find
\begin{equation}
\mapleinline{inert}{2d}{Sum((-1)^n*cos(1/2*Pi*(-2*n^2+4*k*n-k)/(2*k-1)),n = 1 .. 2*k-1) =
0;}{%
\[
{\displaystyle \sum _{n=1}^{2\,k - 1}} \,(-1)^{n}\,\mathrm{cos}(
{\displaystyle \frac {\pi \,( - 2\,n^{2} + 4\,k\,n - k)}{2\,(2\,k
 - 1)}} )=0
\]
}
\label{B18}
\end{equation}

which can alternatively be rewritten as
\begin{equation}
\mapleinline{inert}{2d}{Sum((-1)^n*sin(Pi*n^2/(4*k+1)),n = -2*k .. 2*k) =
Sum((-1)^n*cos(Pi*n^2/(4*k-1)),n = 1-2*k .. 2*k-1);}{%
\[
{\displaystyle \sum _{n= - 2\,k}^{2\,k}} \,(-1)^{n}\,\mathrm{sin}
({\displaystyle \frac {\pi \,n^{2}}{4\,k + 1}} )={\displaystyle 
\sum _{n=1 - 2\,k}^{2\,k - 1}} \,(-1)^{n}\,\mathrm{cos}(
{\displaystyle \frac {\pi \,n^{2}}{4\,k - 1}} )=0\,,
\]
}
\end{equation}

and
\begin{equation}
\mapleinline{inert}{2d}{Sum((-1)^n*sin(1/2*Pi*(-2*n^2+4*k*n-k)/(2*k-1)),n = 1 .. 2*k-1) =
-sqrt(2*k-1);}{%
\[
{\displaystyle \sum _{n=1}^{2\,k - 1}} \,(-1)^{n}\,\mathrm{sin}(
{\displaystyle \frac {\pi \,( - 2\,n^{2} + 4\,k\,n - k)}{2\,(2\,k
 - 1)}} )= - \sqrt{2\,k - 1}\,.
\]
}
\label{B19}
\end{equation}

which is equivalent to
\begin{equation}
\mapleinline{inert}{2d}{L19a := Sum((-1)^n*cos(Pi*n^2/(4*k+1)),n = -2*k .. 2*k) =
(-1)^k*(4*k+1)^(1/2);}{%
\[
{\displaystyle \sum _{n= - 2\,k}^{2\,k}} \,(-1)
^{n}\,\mathrm{cos}({\displaystyle \frac {\pi \,n^{2}}{4\,k + 1}} 
)=(-1)^{k}\,\sqrt{4\,k + 1}
\]
}
\label{L19a}
\end{equation}
or, in a different form,
\begin{equation}
\mapleinline{inert}{2d}{L19b := Sum((-1)^n*sin(Pi*n^2/(4*k-1)),n = 1-2*k .. 2*k-1) =
(-1)^k*(4*k-1)^(1/2);}{%
\[
{\displaystyle \sum _{n=1- 2\,k}^{2\,k - 1}} \,
(-1)^{n}\,\mathrm{sin}({\displaystyle \frac {\pi \,n^{2}}{4\,k - 
1}} )=(-1)^{k}\,\sqrt{4\,k - 1}\,.
\]
}
\label{L19b}
\end{equation}
Similarly, \eqref{S1a=C1a} can be rewritten as

\begin{equation}
\mapleinline{inert}{2d}{Sum(sin(1/4*Pi*(2*n+1)^2/(4*k+1)),n = -2*k .. 2*k) =
(Sum(cos(1/4*Pi*(2*n+1)^2/(4*k+1)),n = -2*k .. 2*k) =
1/2*(8*k+2)^(1/2));}{%
\[
{\displaystyle \sum _{n= - 2\,k}^{2\,k}} \,\mathrm{sin}(
{\displaystyle \frac {\pi \,(2\,n + 1)^{2}}{4\,(4\,k + 1)}} )=
 {\displaystyle \sum _{n= - 2\,k}^{2\,k}} \,\mathrm{
cos}({\displaystyle \frac {\pi \,(2\,n + 1)^{2}}{4\,(4\,k + 1)}} 
)={\displaystyle \frac {\sqrt{8\,k + 2}}{2}}\,.  
\]
}
\label{S5_alt}
\end{equation}
 All of these are reminiscent of special values of a truncated Jacobi theta function
\begin{equation}
\mapleinline{inert}{2d}{theta[T](z,tau) := Sum(exp(i*Pi*(n^2*tau+2*n*z)),n = 0 .. L) =
Sum(cos(Pi*(n^2*tau+2*n*z))+i*sin(Pi*(n^2*tau+2*n*z)),n = 0 .. L);}{%
\[
{\theta _{L}}(z, \,\tau ) \equiv {\displaystyle \sum _{n=-L}^{L}} \,e
^{(i\,\pi \,(n^{2}\,\tau  + 2\,n\,z))}={\displaystyle \sum _{n=-L}
^{L}} \,(\mathrm{cos}(\pi \,(n^{2}\,\tau  + 2\,n\,z)) + i\,
\mathrm{sin}(\pi \,(n^{2}\,\tau  + 2\,n\,z)))\,,
\]
}
\label{theta}
\end{equation}
also recognized as a generalized Gauss quadratic sum \cite[Theorem 9.16] {Apostol}. 
\newline

\subsection{Odd Indexed Quadratic Sums}

Break \eqref{A1R2e} and \eqref{A1R2o} into their even and odd parts to find
\begin{equation}
\mapleinline{inert}{2d}{Sum(sin(Pi*n^2/k),n = 0 .. k)-Sum(sin(1/4*Pi*(2*n-1)^2/k),n = 1 .. k)
= 1/2*(-1)^k*2^(1/2)*k^(1/2);}{%
\[
 {\displaystyle \sum _{n=0}^{k}} \,\mathrm{sin}(
{\displaystyle \frac {\pi \,n^{2}}{k}} )  -  {\displaystyle \sum _{n=1}^{k}} \,\mathrm{sin}(
{\displaystyle \frac {\pi \,(2\,n - 1)^{2}}{4\,k}} )
= {\displaystyle \frac {(-1)^{k}\,\sqrt{2\,k}}{2}} 
\]
}
\label{SinSplit}
\end{equation}
and
\begin{equation}
\mapleinline{inert}{2d}{Sum(cos(Pi*n^2/k),n = 1 .. k)-Sum(cos(1/4*Pi*(2*n-1)^2/k),n = 1 .. k)
= -1/2+1/2*(-1)^k*(1+2^(1/2)*k^(1/2));}{%
\[
{\displaystyle \sum _{n=1}^{k}} \,\mathrm{cos}(
{\displaystyle \frac {\pi \,n^{2}}{k}} )  - {\displaystyle \sum _{n=1}^{k}} \,\mathrm{cos}(
{\displaystyle \frac {\pi \,(2\,n - 1)^{2}}{4\,k}} ) 
= - {\displaystyle \frac {1}{2}}  + {\displaystyle \frac {(-1)^{k
}\,(1 + \sqrt{2\,k})}{2}} 
\]
}
\label{CosSplit}
\end{equation}

Employing the method outlined in Appendix A, and using \eqref{S6p} (see next section below) we find
\begin{equation}
\mapleinline{inert}{2d}{Sum(sin(Pi*n^2/k),n = 0 .. k) = 1/4*2^(1/2)*k^(1/2)*(1+(-1)^k);}{%
\[
{\displaystyle \sum _{n=0}^{k}} \,\mathrm{sin}({\displaystyle 
\frac {\pi \,n^{2}}{k}} )={\displaystyle \frac {\sqrt{2\,k
}\,(1 + (-1)^{k})}{4}} 
\]
}
\label{SinLemma}
\end{equation}
and

\begin{equation}
{\displaystyle \sum _{n=1}^{k}} \,\mathrm{cos}({\displaystyle 
\frac {\pi \,n^{2}}{k}} )= \begin{cases}
   \sqrt{2\,k}/2 & \text{if $k$ is even} \\
   -1       & \text{if $k$ is odd.}
  \end{cases}
\label{CosLemma}
\end{equation}

From \eqref{SinLemma} and \eqref{CosLemma}, the odd indexed quadratic sums appearing in \eqref{SinSplit} and \eqref{CosSplit} eventually give 
\begin{equation}
\mapleinline{inert}{2d}{Sum(sin(1/4*Pi*(2*n-1)^2/k),n = 1 .. k) =
-1/4*2^(1/2)*k^(1/2)*((-1)^k-1);}{%
\[
{\displaystyle \sum _{n=1}^{k}} \,\mathrm{sin}({\displaystyle 
\frac {\pi \,(2\,n - 1)^{2}}{4\,k}} )={\displaystyle \sum _{n=1}^{k}} \,\mathrm{cos}({\displaystyle 
\frac {\pi \,(2\,n - 1)^{2}}{4\,k}} )= - {\displaystyle \frac {
\sqrt{2\,k}\,((-1)^{k} - 1)}{4}}\,, 
\]
\label{OddSinSum}
}
\end{equation}
a generalization of \eqref{S1a=C1a}. This is discussed further in Section 5.

\section{Other Results}
\subsection{Generalizations}
Along these lines, we also note that some of the above may be generalized by extending the upper limit of summation. If $p$ is a positive integer, it may be shown (by induction and removing all integral multiples of $\pi$ from the arguments of the corresponding trigonometric functions - see the Appendix for a sample proof) that \eqref{A1R2e} generalizes to
\begin{equation}
\mapleinline{inert}{2d}{Sum((-1)^n*sin(1/4*Pi*n^2/k),n = 1 .. 4*p*k) =
p*(-1)^k*2^(1/2)*k^(1/2);}{%
\[
{\displaystyle \sum _{n=1}^{4\,k\,p}} \,(-1)^{n}\,\mathrm{sin}(
{\displaystyle \frac {\pi \,n^{2}}{4\,k}} )=p\,(-1)^{k}
\,\sqrt{2\,k}
\]
}
\label{S1pX}
\end{equation}

\eqref{Ans3b} generalizes to

\begin{equation}
\mapleinline{inert}{2d}{Sum((-1)^n*sin(1/2*Pi*n^2/(2*k-1)),n = 1 .. 2*p*(2*k-1)) =
(-1)^k*p*(2*k-1)^(1/2);}{%
\[
{\displaystyle \sum _{n=1}^{2\,(2\,k - 1)\,p}} \,(-1)^{n}\,
\mathrm{sin}({\displaystyle \frac {\pi \,n^{2}}{2\,(2\,k - 1)}} )
=(-1)^{k}\,p\,\sqrt{2\,k - 1}
\]
}
\label{S2pX}
\end{equation}

\eqref{S1a} generalizes to
\begin{equation}
\mapleinline{inert}{2d}{Sum((-1)^n*sin(Pi*n^2/(4*k-1)),n = 1 .. p*(4*k-1)-2*k) =
(-1)^k*(p-1/2)*(4*k-1)^(1/2);}{%
\[
{\displaystyle \sum _{n=1}^{(4\,k - 1)\,p}} \,(-1)^{n}\,
\mathrm{sin}({\displaystyle \frac {\pi \,n^{2}}{4\,k - 1}} )=(-1)
^{k}\,p \,\sqrt{4\,k - 1}\,,
\]
}
\label{S3p}
\end{equation}

and \eqref{S4a} generalizes to
\begin{equation}
\mapleinline{inert}{2d}{Sum((-1)^n*sin(Pi*n^2/(4*k-3)),n = 1 .. p*(4*k-3)-2*k+1) = 0;}{%
\[
{\displaystyle \sum _{n=1}^{(4\,k - 3)\,p}} \,(-1)^{n}
\,\mathrm{sin}({\displaystyle \frac {\pi \,n^{2}}{4\,k - 3}} )=0\,.
\]
}
\label{S4p}
\end{equation}

Similar properties also hold for the alternating cosine sums.
\eqref{A1R2o} generalizes to

\begin{equation}
\mapleinline{inert}{2d}{Sum((-1)^n*cos(1/2*Pi*n^2/(2*k-1)),n = 1 .. 2*p*(2*k-1)) =
-(-1)^k*p*(2*k-1)^(1/2);}{%
\[
{\displaystyle \sum _{n=1}^{2\,(2\,k - 1)\,p}} \,(-1)^{n}\,
\mathrm{cos}({\displaystyle \frac {\pi \,n^{2}}{2\,(2\,k - 1)}} )
= - (-1)^{k}\,p\,\sqrt{2\,k - 1}
\]
}
\label{B3aX}
\end{equation}
 \eqref{Ans1b} generalizes to

\begin{equation}
\mapleinline{inert}{2d}{B4p8x := Sum((-1)^n*cos(1/4*Pi*n^2/k),n = 1 .. 4*p*k) =
p*(-1)^k*2^(1/2)*k^(1/2);}{%
\[
{\displaystyle \sum _{n=1}^{4\,k\,p}} \,(-1)^{n
}\,\mathrm{cos}({\displaystyle \frac {\pi \,n^{2}}{4\,k}} )=p\,(
-1)^{k}\,\sqrt{2\,k}
\]
}
\label{B4aX}
\end{equation}

the variation \eqref{C1a} generalizes to
\begin{equation}
\mapleinline{inert}{2d}{Sum((-1)^n*cos(Pi*n^2/(4*k-1)), n = 1 .. 4*p*k-p) = 0}{\[\displaystyle \sum _{n=1}^{(4\,k-1)\,p} \left( -1 \right) ^{n}\cos \left( {\frac {\pi \,{n}^{2}}{4\,k-1}} \right) =0\]}\, ,
\label{S1px_var}
\end{equation}

and \eqref{S3a} generalizes to
\begin{equation}
\mapleinline{inert}{2d}{B4p11P := Sum((-1)^n*cos(Pi*n^2/(4*k-3)),n = 1 .. p*(4*k-3)) =
-(-1)^k*p*(4*k-3)^(1/2);}{%
\[
{\displaystyle \sum _{n=1}^{(4\,k - 3)\,p}} \,
(-1)^{n}\,\mathrm{cos}({\displaystyle \frac {\pi \,n^{2}}{4\,k - 
3}} )= - (-1)^{k}\,p\,\sqrt{4\,k - 3}\,.
\]
}
\label{C3p}
\end{equation}
\newline

\subsection{Other Miscellaneous results}
Ancillary to proofs of the above, we obtained a number of miscellaneous results, all of which can be proven by using the methods outlined in the previous section with the help of the Appendix. To prove \eqref{S2p} for $p$ odd, the following companion to \eqref{A1R2o} is needed 
\begin{equation}
\mapleinline{inert}{2d}{Sum(cos(1/2*Pi*n^2/(2*k-1)),n = 1 .. 2*k-1) =
1/2*(2*k-1)^(1/2)-1/2;}{%
\[
{\displaystyle \sum _{n=1}^{2\,k - 1}} \,\mathrm{cos}(
{\displaystyle \frac {\pi \,n^{2}}{2\,(2\,k - 1)}} )=
{\displaystyle \frac {\sqrt{2\,k - 1}}{2}}  - {\displaystyle 
\frac {1}{2}} \, ,
\]
}
\label{S2p_C}
\end{equation}
and the following, needed for the proof of \eqref{B4a} is a companion to \eqref{Ans1b}
\begin{equation}
\mapleinline{inert}{2d}{Sum(cos(1/4*Pi*n^2/k),n = 1 .. 2*k) =
half*(2^(1/2)*k^(1/2)+(-1)^k-1);}{%
\[
{\displaystyle \sum _{n=1}^{2\,k}} \,\mathrm{cos}({\displaystyle 
\frac {\pi \,n^{2}}{4\,k}} )={\displaystyle \frac{\sqrt{2\,k}
 + (-1)^{k} - 1}{2}}\, .
\]
}
\label{lemma2}
\end{equation}

Further, based on \eqref{A1R2o} and \eqref{Hs1} with $k \rightarrow 4(2k-1)$ followed by a lengthy series of reductions of the upper limit, we have
\begin{equation}
\mapleinline{inert}{2d}{Sum(sin(1/2*Pi*n^2/(2*k-1)),n = 1 .. 2*k-2) =
1/2*(2*k-1)^(1/2)+1/2*(-1)^k;}{%
\[
{\displaystyle \sum _{n=1}^{2\,k - 1}} \,\mathrm{sin}(
{\displaystyle \frac {\pi \,n^{2}}{2\,(2\,k - 1)}} )=
{\displaystyle \frac {\sqrt{2\,k - 1}}{2}}  - {\displaystyle 
\frac {(-1)^{k}}{2}}\,. 
\]
}
\label{Bc3_lemma}
\end{equation}

Further, \eqref{S1a=C1a} generalizes to
\begin{equation}
\mapleinline{inert}{2d}{Sum(sin(1/4*Pi*(2*n-1)^2/(2*k-1)),n = 1 .. p*(2*k-1)) =
1/2*p*(4*k-2)^(1/2);}{%
\[
{\displaystyle \sum _{n=1}^{p\,(2\,k - 1)}} \,\mathrm{sin}(
{\displaystyle \frac {\pi \,(2\,n - 1)^{2}}{4\,(2\,k - 1)}} )=
{\displaystyle \frac {p\,\sqrt{4\,k - 2}}{2}}, 
\]
}
\label{S5p}
\end{equation}
and \eqref{Hs1} generalizes to
\begin{equation}
\mapleinline{inert}{2d}{Sum(sin(2*Pi*n^2/k),n = 1 .. p*k-1) =
1/2*p*k^(1/2)*(1+cos(1/2*k*Pi)-sin(1/2*k*Pi));}{%
\[
{\displaystyle \sum _{n=1}^{p\,k }} \,\mathrm{sin}(
{\displaystyle \frac {2\,\pi \,n^{2}}{k}} )={\displaystyle 
\frac {1}{2}} \,p\,\sqrt{k}\,(1 + \mathrm{cos}({\displaystyle 
\frac {k\,\pi }{2}} ) - \mathrm{sin}({\displaystyle \frac {k\,\pi
 }{2}} ))\,.
\]
}
\label{S6p}
\end{equation}
\eqref{A1R2e} generalizes to 
\begin{equation}
\mapleinline{inert}{2d}{Sum((-1)^n*sin(1/4*Pi*n^2/k),n = 1 .. 2*p*k-1) =
1/2*p*(-1)^k*2^(1/2)*k^(1/2);}{%
\[
{\displaystyle \sum _{n=1}^{2\,p\,k}} \,(-1)^{n}\,\mathrm{sin
}({\displaystyle \frac {\pi \,n^{2}}{4\,k}} )={\displaystyle 
 {p\,(-1)^{k}\,\sqrt{k/2}}{}}\, , 
\]
}
\label{S1p}
\end{equation}
\eqref{A1R2o} generalizes to

\begin{equation}
\mapleinline{inert}{2d}{Sum((-1)^n*cos(1/2*Pi*n^2/(2*k-1)),n = 1 .. p*(2*k-1)) =
-(1-(-1)^p)/four-1/2*(-1)^k*p*(2*k-1)^(1/2);}{%
\[
{\displaystyle \sum _{n=1}^{p\,(2\,k - 1)}} \,(-1)^{n}\,\mathrm{
cos}({\displaystyle \frac {\pi \,n^{2}}{2\,(2\,k - 1)}} )= - 
{\displaystyle \frac {1 - (-1)^{p}}{{4}}}  - 
{\displaystyle \frac {(-1)^{k}\,p\,\sqrt{2\,k - 1}}{2}}\,, 
\]
}
\label{B3a}
\end{equation}

\eqref{Ans1b} generalizes to
\begin{equation}
\mapleinline{inert}{2d}{Sum((-1)^n*cos(1/4*Pi*n^2/k),n = 1 .. 2*p*k) =
1/2*p*(-1)^k*2^(1/2)*k^(1/2)-(1-(-1)^(p*k))/two;}{%
\[
{\displaystyle \sum _{n=1}^{2\,p\,k}} \,(-1)^{n}\,\mathrm{cos}(
{\displaystyle \frac {\pi \,n^{2}}{4\,k}} )={\displaystyle 
\frac {p\,(-1)^{k}\,\sqrt{2\,k}}{2}}  - {\displaystyle 
\frac {1 - (-1)^{p\,k}}{{2}}} \,,
\]
}
\label{B4a}
\end{equation}

\eqref{Ans3b} generalizes to
\begin{equation}
\mapleinline{inert}{2d}{Sum((-1)^n*sin(1/2*Pi*n^2/(2*k-1)),n = 1 .. p*(2*k-1)) =
1/4*(-1)^k*(1+2*p*(2*k-1)^(1/2)-(-1)^p);}{%
\[
{\displaystyle \sum _{n=1}^{p\,(2\,k - 1)}} \,(-1)^{n}\,\mathrm{
sin}({\displaystyle \frac {\pi \,n^{2}}{2\,(2\,k - 1)}} )=
{\displaystyle \frac {(-1)^{k}\,(1 + 2\,p\,\sqrt{2\,k - 1} - (-1)
^{p})}{4}}\,. 
\]
}
\label{S2p}
\end{equation}

Along with the above, we have also found

\begin{equation}
\mapleinline{inert}{2d}{Sum((-1)^n*cos(Pi*n^2/(4*k-1)),n = 1 .. p*(4*k-1)-2*k) = -1/2;}{%
\[
{\displaystyle \sum _{n=1}^{p\,(4\,k - 1) - 2\,k}} \,(-1)^{n}\,
\mathrm{cos}({\displaystyle \frac {\pi \,n^{2}}{4\,k - 1}} )=
{\displaystyle -\frac {1}{2}}\,, 
\]
}
\label{S1px}
\end{equation}

\begin{equation} 
\mapleinline{inert}{2d}{Sum((-1)^n*cos(Pi*n^2/(4*k-3)),n = 1 .. p*(4*k-3)-2*k+1) =
-1/2-(-1)^k*(4*k-3)^(1/2)*(p-1/2);}{%
\[
{\displaystyle \sum _{n=1}^{p\,(4\,k - 3) - 2\,k + 1}} \,(-1)^{n}
\,\mathrm{cos}({\displaystyle \frac {\pi \,n^{2}}{4\,k - 3}} )=
 - {\displaystyle \frac {1}{2}}  - (-1)^{k}\,\sqrt{4\,k - 3}\,(p
 - {\displaystyle \frac {1}{2}} )\,,
\]
\label{B12b}
}
\end{equation}
and
\begin{equation}
\mapleinline{inert}{2d}{B28g := Sum((-1)^n*sin(Pi*n^2/(4*k-1)),n = 1 .. p*(4*k-1)-2*k) =
(-1)^k*(p-1/2)*(4*k-1)^(1/2);}{%
\[
{\displaystyle \sum _{n=1}^{p\,(4\,k - 1) - 2\,k
}} \,(-1)^{n}\,\mathrm{sin}({\displaystyle \frac {\pi \,n^{2}}{4
\,k - 1}} )=(-1)^{k}\,(p - {\displaystyle \frac {1}{2}} )\,\sqrt{
4\,k - 1} \,.
\]
}
\label{B28g}
\end{equation}

All of the above may be proven using Appendix A as a template.
\section{Connection with Gauss Quadratic Sums}

As usual, denote the ``Extended Gauss Quadratic Sum'' for $j,k,m,p \in \mathbb{N}$ by

\begin{equation}
\mapleinline{inert}{2d}{G(`j;m`) = [Sum(exp(2*I*Pi*j*n^2/k),n = 1 .. k) =
Sum(cos(2*Pi*j*n^2/m),n = 1 .. m)+Sum(sin(2*Pi*j*n^2/m),n = 1 ..
m)*I];}{%
\[
\mathit{G_{p}}(\mathit{j;k;\theta})= {\displaystyle \sum _{n=1}^{
k\,p}} \,\exp{\displaystyle (\frac {2\,i\,\pi \,j\,n^{2}}{k}+2\,\pi\,i\,\theta \, n)}\,,
\]
}
\label{GaussSum}
\end{equation}
and write the (classical) ``Gauss Quadratic Sums'' as 
\begin{equation}
G_{1}(j;k;0)= 
{\displaystyle \sum _{n=1}^{k}} \,\mathrm{cos}({\displaystyle 
\frac {2\,\pi \,j\,n^{2}}{k}} ) +  
{\displaystyle i\,\sum _{n=1}^{k}} \,\mathrm{sin}({\displaystyle 
\frac {2\,\pi \,j\,n^{2}}{k}} )\,. 
\end{equation}

Consistent with this notation, we define the ``Extended Alternating Gauss Quadratic Sum'' by

\begin{equation}
\mapleinline{inert}{2d}{G^{A}_{p}(`j;m`) = [Sum((-1)^n*exp(2*I*Pi*j*n^2/m),n = 1 .. m) =
Sum((-1)^n*cos(2*Pi*j*n^2/m),n = 1 ..
m)+Sum((-1)^n*sin(2*Pi*j*n^2/m),n = 1 .. m)*I];}{%
\[
{G^{A}_{p}}(\mathit{j;k})\equiv {\displaystyle \sum _{n=1}^{p\,k}
} \,(-1)^{n}\,\exp{\displaystyle(\frac {i\,\pi \,j\,n^{2}}{k})}=
G_{p/2}(j;2k;1/2)
%
\]
}
\label{AltGaussSum}
\end{equation}

Along with Gauss' classical result, equivalent to \eqref{Hs1} and \eqref{Hc1}, (one of whose proofs involves contour integrals similar to the method used to obtain our results quoted in Section 2)

\begin{equation}
\mapleinline{inert}{2d}{G("1;m") = 1/2*m^(1/2)*(1+i)*(1+exp(-1/2*i*Pi*m));}{%
\[
\mathit{G_{1}}({1;k;0})={\displaystyle \frac {1}{2}} \,\sqrt{k
}\,(1 + i)\,(1 + \exp{\displaystyle( - {i\,\pi \,k/2}{})})\, ,
\]
}
\label{G(1;k)}
\end{equation}

the entire content of Section 4.1 can be newly summarized as a special case of \eqref{AltGaussSum} - in terms of four roots of unity - by 
\begin{equation}
\mapleinline{inert}{2d}{G[A]("1/2;m") = m^(1/2)*exp(1/4*I*Pi*(1-m));}{%
\[
{G^{A}_{p}}({1;k})={\displaystyle \sum _{n=1}^{p\,k}
} \,(-1)^{n}\,\exp{\displaystyle(\frac {i\,\pi \,n^{2}}{k})}=p\,\sqrt{k}\,\exp{(i\pi \,(1/4+\lfloor{k/4}\rfloor - \mod (k,4)/4))}
\]\, .
}
\label{GA(1/2;k}
\end{equation}

As an extension to \eqref{G(1;k)}, with reference to \eqref{S1a=C1a},\eqref{OddSinSum} and \eqref{S5p}, we also obtain the following Generalized Odd-indexed Quadratic Gauss sum
\begin{equation}
\mapleinline{inert}{2d}{Ge1 := sum(exp(1/4*I*Pi*(2*n-1)^2/m),n = 1 .. p*m) =
(1+I)*p/four*2^(1/2)*m^(1/2)*(1-(-1)^m);}{%
\[
{\displaystyle \sum _{n=1}^{p\,k}} \,\exp({\displaystyle \frac {i\,\pi \,(2\,n - 1)^{2}}{4\,k}})={\displaystyle  {(1 + i)\,
p\,\sqrt{2\,k}\,(1 - (-1)^{k})}} /4 \, ,
\]
}
\label{Ge1}
\end{equation}
which could be compared to the known result \cite{Waterhouse}
\begin{equation}
\mapleinline{inert}{2d}{Sum(exp(2*I*Pi*(n-k)^2/four/r),n = 0 .. 2*r-1) = (two*r*I)^(1/2);}{%
\[
{\displaystyle \sum _{n=0}^{2\,k - 1}} \,\exp{\displaystyle (\frac {i\,\pi \,(
n - m)^{2}}{2\,k})}=\sqrt{2\,k\,i}\,.
\]
}
\label{Water1}
\end{equation}

Write \eqref{Water1} in the form of \eqref{Ge1} by setting $k \rightarrow k/2, m=1/2$ and noting that the result is formally only valid (if it is technically valid at all) for $k$ even suggests that


\begin{equation}
\mapleinline{inert}{2d}{Sum(exp(i*Pi*(n-1/2)^2/k),n = 0 .. k-1) = (i*m)^(1/2);}{%
\[
{\displaystyle \sum _{n=0}^{k - 1}} \,\exp{\displaystyle(\frac {i\,\pi \,(n - 1/2)^{2}}{k})}=\sqrt{i\,k}\, ,
\]
}
\label{Water2}
\end{equation}
which surprisingly reduces to \eqref{Ge1} for {\it{odd}} values of $k$ when $p=1$, after correcting for the difference in the summation limits between the two results. {\emph {Although interesting as a special limiting case of \eqref{Ge1}, \eqref{Water2} is only correct when $k$ is odd}}.\newline

According to our notational definitions, an equivalent form of \eqref{Ge1} is
\begin{equation}
\mapleinline{inert}{2d}{Ge1a := sum(exp(i*Pi*(n^2-n)/m),n = 1 .. p*m) =
1/4*(1+i)*exp(-1/4*i*Pi/m)*p*m^(1/2)*2^(1/2)*(1-(-1)^m);}{%
\[
{G_{p/2}(1;2\,k;-\frac{1}{2\,k}) = \displaystyle \sum _{n=1}^{p\,k}} \,\exp{\displaystyle (\frac {
i\,\pi \,(n^{2} - n)}{k})}={\displaystyle \frac {1}{4}} \,(1 + i)
\,\exp{\displaystyle( - {i\,\pi }/(4\,k){})}\,p\,\sqrt{2\,k}\,(1 - (-1
)^{k})\,,
\]
}
\label{Ge1a}
\end{equation}
the asymptotics of which have recently been revisited \cite{Paris} without reference to the closed form \eqref{Ge1a}. An alternate (canonical) form of \eqref{Ge1} is
\begin{equation}
\mapleinline{inert}{2d}{Eq5p9b := Sum(exp(Pi*(n^2-n)/k*I),n = 0 .. p*(k-1)) =
(1+I)*exp(-I/four/k*Pi)*p*sqrt(2*k)*(1-(-1)^k)/four+1-(-1)^(p^2*k+p)*S
um(exp(Pi*n*(n-1)/k*I),n = 1 .. p);}{%
\maplemultiline{
{\displaystyle \sum _{n=0}^{p\,(k - 1)}} \,\exp{ \displaystyle
(\frac {i\,\pi \,(n^{2} - n)}{k})}= 
{\displaystyle \frac{p}{4} {\,(1 + i)\,\exp{\displaystyle ({ - i\,\pi/4k })}\,\sqrt{2\,k}\,(1 - (-1)^{k})}} \\ +\, 1\,
 - (-1)^{(p^{2}\,k + p)}\, {\displaystyle \sum _{n=1}
^{p}} \,\exp{(\displaystyle \frac {i\,\pi \,n\,(n - 1)}{k})} \,. }
}
\label{canonical}
\end{equation}
We emphasize that \eqref{canonical} is valid for all integers $p,k$ and thereby generalizes the (ought-to-be) known result
\begin{equation}
\mapleinline{inert}{2d}{Sum(exp(Pi*(n^2-n)/k*I),n = 0 .. k-1) =
(1/2+1/2*I)*exp(-1/4*I/k*Pi)*2^(1/2)*k^(1/2);}{%
\[
{\displaystyle \sum _{n=0}^{k - 1}} \,\exp{\displaystyle (\frac {i\,\pi \,(n^{2} - n
)}{k})}={\displaystyle \frac{1}{2}}({\displaystyle 1}  + {\displaystyle } \,i)\,\exp{\displaystyle (\frac {-i\,\pi }{4\,k})}\,
\sqrt{2\,k}
\]
}
\label{LsL=LsR}
\end{equation}
which is only valid for odd values of $k$. Since \eqref{LsL=LsR} appears not to be well-known, we note that it can be simply obtained by applying the generalized Landsberg-Schaar identity \cite[Eq.(2.8)]{ BerEvans}
\begin{equation}
\mapleinline{inert}{2d}{j^(1/2)*Sum(exp(i*Pi*(j*n^2+n*m)/k),n = 0 .. k-1) =
k^(1/2)*exp(1/4*i*Pi*(j*k-m^2)/j/k)*Sum(exp(-i*Pi*(n^2*k+n*m)/j),n = 0
.. j-1);}{%
\[
\sqrt{j}\, {\displaystyle \sum _{n=0}^{k - 1}} \,\exp{\displaystyle (
\frac {i\,\pi \,(j\,n^{2} + n\,m)}{k})}  =\sqrt{k}\,\exp
{\displaystyle (\frac {i\,\pi \,(j\,k - m^{2})}{4\,j\,k})}\, 
{\displaystyle \sum _{n=0}^{j - 1}} \,\exp{\displaystyle( - \frac {i\,\pi \,(n^{
2}\,k + n\,m)}{j})}  
\]
}
\label{LsIdent}
\end{equation}
to the left-hand side of \eqref{LsL=LsR}, with the caveat that \eqref{LsIdent} is only valid when $jk+m$ is even. Utilizing the freedom to choose $p$ leads to further interesting results. For example, let $p \rightarrow p\,k$ in \eqref{canonical} and, with reference to \eqref{Ge1a}, obtain
\begin{equation}
\mapleinline{inert}{2d}{Eq9p5c := Sum(exp(Pi*(n^2-n)/k*I),n = 0 .. p*k*(k-1)) =
1+1/4*p*(1+i)*exp(-1/4*i*Pi/k)*sqrt(2*k)*(k-1)*(1-(-1)^k);}{%
\[
{\displaystyle \sum _{n=0}^{p\,k\,(k - 1)}} \,
\exp{(\displaystyle \frac {i\,\pi \,(n^{2} - n)}{k})}=1 + {\displaystyle \frac {1
}{4}} \,p\,(1 + i)\,\exp{\displaystyle( - \frac {i\,\pi }{4\,k})}\,\sqrt{2\,k}\,
(k - 1)\,(1 - (-1)^{k})\,.
\]
}
\label{var1}
\end{equation}
Alternatively, let $p\rightarrow p\,(k-1)$ to obtain
\begin{equation}
\mapleinline{inert}{2d}{Eq5p9e := Sum(exp(i*Pi*n*(n-1)/k),n = 0 .. p*(k-1)^2) =
1+Sum(exp(i*Pi*n*(n-1)/k),n = 1 ..
p)+1/4*(1+i)*(k-2)*p*sqrt(2*k)*exp(-1/4*i*Pi/k)*(1-(-1)^k);}{%
\maplemultiline{
{\displaystyle \sum _{n=0}^{p\,(k - 1)^{2}}} 
\,\exp{\displaystyle(\frac {i\,\pi \,n\,(n - 1)}{k})}=1 + {\displaystyle \sum _{n=1}^{p}} \,\exp{\displaystyle(\frac {i\,\pi \,n\,(n - 1)
}{k})}  \\
\hspace{2cm}+\,{\displaystyle \frac {1}{4}} \,(1 + i)\,(k
 - 2)\,p\,\sqrt{2\,k}\,\exp{\displaystyle( - \frac {i\,\pi }{4\,k})}\,(1 - (-1)
^{k})\,.} 
}
\label{var2}
\end{equation}

\section{Acknowledgements} We thank N.E. Frankel who drew our attention to \cite{HardyLitt}.

\section{Appendix}

The following proof of \eqref{S3p} is a template for the proof of all results presented in Sections 3.2 and 4. We proceed by induction with the claim

\begin{maplegroup}

\mapleresult
\begin{maplelatex}
\begin{equation}
\mapleinline{inert}{2d}{Sum((-1)^n*sin(Pi*n^2/(4*k-1)),n = 1 .. p*(4*k-1)) =
(-1)^k*p*(4*k-1)^(1/2);}{%
\[
{\displaystyle \sum _{n=1}^{p\,(4\,k - 1)}} \,(-1)^{n}\,\mathrm{
sin}({\displaystyle \frac {\pi \,n^{2}}{4\,k - 1}} )=(-1)^{k}\,p
\,\sqrt{4\,k - 1} .
\]
}
\label{S3pX}
\end{equation}
\end{maplelatex}

\end{maplegroup}

The case $p=1$ is known to be true from \eqref{S1a}. Set $p=2$ in \eqref{S3pX} giving 
\begin{maplegroup}

\mapleresult
\begin{maplelatex}
\begin{equation}
\mapleinline{inert}{2d}{B28b := Sum((-1)^n*sin(Pi*n^2/(4*k-1)),n = 1 .. 8*k-2) =
2*(-1)^k*(4*k-1)^(1/2);}{%
\[
{\displaystyle \sum _{n=1}^{8\,k - 2}} \,(-1)^{n
}\,\mathrm{sin}({\displaystyle \frac {\pi \,n^{2}}{4\,k - 1}} )=2
\,(-1)^{k}\,\sqrt{4\,k - 1}\,.
\]
}
\label{S3pXa}
\end{equation}
\end{maplelatex}

\end{maplegroup}

Splitting the sum into two parts gives
\begin{maplegroup}

\mapleresult
\begin{maplelatex}
\begin{equation}
\mapleinline{inert}{2d}{B28b := Sum((-1)^n*sin(Pi*n^2/(4*k-1)),n = 1 ..
4*k-1)+Sum((-1)^n*sin(Pi*n^2/(4*k-1)),n = 4*k .. 8*k-2) =
2*(-1)^k*(4*k-1)^(1/2);}{%
\[
{\displaystyle \sum _{n=1}^{4\,k - 1
}} \,(-1)^{n}\,\mathrm{sin}({\displaystyle \frac {\pi \,n^{2}}{4
\,k - 1}} ) +  {\displaystyle \sum _{n=4
\,k}^{8\,k - 2}} \,(-1)^{n}\,\mathrm{sin}({\displaystyle \frac {
\pi \,n^{2}}{4\,k - 1}} ) =2\,(-1)^{k}\,\sqrt{4\,k - 
1}
\]
}
\label{S3pXb}
\end{equation}
\end{maplelatex}

\end{maplegroup}
and shifting the summation indices of the second sum, yields
\begin{maplegroup}

\mapleresult
\begin{maplelatex}
\begin{equation}
\mapleinline{inert}{2d}{B28b := Sum((-1)^n*sin(Pi*n^2/(4*k-1)),n = 1 ..
4*k-1)+Sum((-1)^(n+4*k-1)*sin(Pi*(n+4*k-1)^2/(4*k-1)),n = 1 .. 4*k-1)
= 2*(-1)^k*(4*k-1)^(1/2);}{%
\maplemultiline{
{\displaystyle \sum _{n=1}^{4\,k - 1
}} \,(-1)^{n}\,\mathrm{sin}({\displaystyle \frac {\pi \,n^{2}}{4
\,k - 1}} ) + {\displaystyle \sum _{n=1}
^{4\,k - 1}} \,(-1)^{(n + 4\,k - 1)}\,\mathrm{sin}(
{\displaystyle \frac {\pi \,(n + 4\,k - 1)^{2}}{4\,k - 1}} ) = 
2\,(-1)^{k}\,\sqrt{4\,k - 1}\,. }
}
\label{S3pXc}
\end{equation}
\end{maplelatex}

\end{maplegroup}

Identify the first sum via \eqref{S1a}, square the argument of the sine term and remove common factors in the second sum to obtain
\begin{maplegroup}

\mapleresult
\begin{maplelatex}
\begin{equation}
\mapleinline{inert}{2d}{B28b :=
(-1)^k*(4*k-1)^(1/2)+Sum((-1)^(n-1)*sin(Pi*(n^2/(4*k-1)+2*n+4*k-1)),n
= 1 .. 4*k-1) = 2*(-1)^k*(4*k-1)^(1/2);}{%
\maplemultiline{
(-1)^{k}\,\sqrt{4\,k - 1} +  {\displaystyle \sum _{n=1
}^{4\,k - 1}} \,(-1)^{(n - 1)}\,\mathrm{sin}(\pi \,(
{\displaystyle \frac {n^{2}}{4\,k - 1}}  + 2\,n + 4\,k - 1)) 
=2\,(-1)^{k}\,\sqrt{4\,k - 1} }
}
\label{S3pxd}
\end{equation}
\end{maplelatex}

\end{maplegroup}

and after removing integral multiples of $\pi$ find
\begin{maplegroup}

\mapleresult
\begin{maplelatex}
\begin{equation}
\mapleinline{inert}{2d}{B28b := (-1)^k*(4*k-1)^(1/2)+Sum((-1)^n*sin(Pi*n^2/(4*k-1)),n = 1 ..
4*k-1) = 2*(-1)^k*(4*k-1)^(1/2);}{%
\[
(-1)^{k}\,\sqrt{4\,k - 1} +  
{\displaystyle \sum _{n=1}^{4\,k - 1}} \,(-1)^{n}\,\mathrm{sin}(
{\displaystyle \frac {\pi \,n^{2}}{4\,k - 1}} ) =2\,(
-1)^{k}\,\sqrt{4\,k - 1}
\]
}
\label{S3pxe}
\end{equation}
\end{maplelatex}
\end{maplegroup}

yielding the requisite identity because of \eqref{S1a}. Therefore, \eqref{S3pX} is true for $p=1$ and $p=2$. \newline

Assuming \eqref{S3pX} for any $p$, let $p \rightarrow p+1$, giving
\begin{equation}
\mapleinline{inert}{2d}{Sum((-1)^n*sin(Pi*n^2/(4*k-1)),n = 1 .. (p+1)*(4*k-1)) =
(-1)^k*(p+1)*(4*k-1)^(1/2);}{%
\[
{\displaystyle \sum _{n=1}^{(p + 1)\,(4\,k - 1)}} \,(-1)^{n}\,
\mathrm{sin}({\displaystyle \frac {\pi \,n^{2}}{4\,k - 1}} )=(-1)
^{k}\,(p + 1)\,\sqrt{4\,k - 1}\,;
\]
}
\label{S3pxf}
\end{equation}
split the summation into two parts as follows
\begin{equation}
\mapleinline{inert}{2d}{Sum((-1)^n*sin(Pi*n^2/(4*k-1)),n = 1 ..
p*(4*k-1))+Sum((-1)^n*sin(Pi*n^2/(4*k-1)),n = p*(4*k-1) ..
(p+1)*(4*k-1)) = (-1)^k*(p+1)*(4*k-1)^(1/2);}{%
\[
{\displaystyle \sum _{n=1}^{p\,(4\,k - 1)}} \,(-1)^{n
}\,\mathrm{sin}({\displaystyle \frac {\pi \,n^{2}}{4\,k - 1}} )
 + {\displaystyle \sum _{n=p\,(4\,k - 1)
}^{(p + 1)\,(4\,k - 1)}} \,(-1)^{n}\,\mathrm{sin}({\displaystyle 
\frac {\pi \,n^{2}}{4\,k - 1}} ) =(-1)^{k}\,(p + 1)\,
\sqrt{4\,k - 1}\,,
\]
}
\label{S3pxg}
\end{equation}
apply \eqref{S3pX} to the first term, shift the summation index of the second, expand and simplify as in the transition from \eqref{S3pXc} to \eqref{S3pxd} to obtain
\begin{equation}
\mapleinline{inert}{2d}{(-1)^k*p*(4*k-1)^(1/2)-(-1)^p*Sum((-1)^(n+1)*sin(Pi*(n^2/(4*k-1)+2*n*
p+p^2)),n = 0 .. 4*k-1) = (-1)^k*(p+1)*(4*k-1)^(1/2);}{%
\[
(-1)^{k}\,p\,\sqrt{4\,k - 1} - (-1)^{p}\, 
{\displaystyle \sum _{n=0}^{4\,k - 1}} \,(-1)^{(n + 1)}\,\mathrm{
sin}(\pi \,({\displaystyle \frac {n^{2}}{4\,k - 1}}  + 2\,n\,p + 
p^{2})) =(-1)^{k}\,(p + 1)\,\sqrt{4\,k - 1}\,.
\]
}
\label{Sp3xh}
\end{equation}
Now, remove integral multiples of $\pi$ from the argument of the sine function, leaving 
\begin{equation}
\mapleinline{inert}{2d}{(-1)^k*p*(4*k-1)^(1/2)+(-1)^(p+p^2)*Sum((-1)^n*sin(Pi*n^2/(4*k-1)),n
= 1 .. 4*k-1) = (-1)^k*(p+1)*(4*k-1)^(1/2);}{%
\[
(-1)^{k}\,p\,\sqrt{4\,k - 1} + (-1)^{(p + p^{2})}\, 
{\displaystyle \sum _{n=1}^{4\,k - 1}} \,(-1)^{n}\,\mathrm{sin}(
{\displaystyle \frac {\pi \,n^{2}}{4\,k - 1}} ) =(-1)
^{k}\,(p + 1)\,\sqrt{4\,k - 1}
\]
}
\label{Sp3xi}
\end{equation}
which reduces to an identity because of \eqref{S3pX} and the fact that  \begin{equation}\nonumber
(-1)^{p+p^2}=1 
\end{equation}
for all integers $p$. Hence by the usual logic of induction \eqref{S3pX} is true.  QED.
Note: In some proofs, it is expedient to reverse the order of summation.

\end{flushleft}
\end{document}